\documentclass{ifacconf}
\usepackage{lipsum}
\usepackage[utf8]{inputenc}
\usepackage{setspace}
\usepackage{natbib}
\usepackage{natbib}
\usepackage{graphicx}        % standard LaTeX graphics tool
\usepackage{subfigure}
\usepackage{amssymb,url}
\usepackage{amsmath}
\newtheorem{prop1}{\bf Proposition}
\newtheorem{lemma1}{\bf Lemma}

\begin{document}
\begin{frontmatter}

\title{Optimal Control of  Joint  Multi-Virus Infection and Information Spreading}

\author[First]{Vladislav~Taynitskiy}
\author[First]{Elena~Gubar}
\author[Third]{Denis~Fedyanin}
\author[Third]{Ilya~Petrov }
\author[Fifth]{Quanyan~Zhu}

\address[First]{St. Petersburg State University, Faculty of Applied
Mathematics and Control Processes, Saint-Petersburg, Russia. (E-mail: tainitsky@gmail.com, e.gubar@spbu.ru) }
%\address[Second]{Department of Electrical and Computer Engineering, Tandon School of Engineering, New York University, Brooklyn, NY, USA, 11201. (E-mail: quanyan.zhu@nyu.edu)}
\address[Third]{V. A. Trapeznikov Institute of Control Sciences of RAS,
Moscow, Russia. (E-mail: dfedyanin@inbox.ru, zyxzy@protonmail.ch)}
%\address[Third]{V. A. Trapeznikov Institute of Control Sciences of RAS,
%Moscow, Russia. (E-mail:
%zyxzy@protonmail.ch)}
\address[Fifth]{Department of Electrical and Computer Engineering,  New York University, USA. (E-mail: quanyan.zhu@nyu.edu)}
\begin{abstract}
Nowadays, epidemic models provide an appropriate tool for describing the propagation of biological viruses in human or animal populations, or rumours and other kinds of information in social networks and malware in both computer and ad hoc networks. Commonly, there are exist multiple types of malware infecting a network of computing devices, or different messages can spread over the social network. Information spreading and virus propagation are interdependent processes. To capture such independencies, we integrate two epidemic models into one holistic framework, known as the modified Susceptible-Warned-Infected-Recovered-Susceptible (SWIRS) model. The first epidemic model describes the information spreading regarding the risk of malware attacks and possible preventive procedures. The second one describes the propagation of multiple viruses over the network of devices.  To minimize the impact of the virus spreading and improve the protection of the networks, we consider an optimal control problem with two types of control strategies: information spreading among healthy nodes and the treatment of infected nodes. We obtain the structure of optimal control strategies and study the condition of epidemic outbreaks. The main results are extended to the case of the network of two connected clusters. Numerical examples are used to corroborate the theoretical findings.
\end{abstract}


\end{frontmatter}

\section{Introduction}\label{introduction}

Recent advances in information technologies have witnessed an exponential growth in the number of devices connected to the Internet and the rapid expansion of the use of social  networks.  The proliferation of devices creates opportunities to spread information more conveniently but has also created a large attack surface for the malware to exploit existing vulnerabilities of the devices and spread malicious codes over the Internet. The channels of malware spreading nowadays are not just limited to computer networks but also include mobile networks and online social networks. Moreover, the wide applications of networks generate an increasing amount of security threats. Computer virus or malware spread and attack a large number of nodes as the network connectivity increases. It can disrupt computer functionalities, collect sensitive confidential information, and gain illegal access to private computer networks at a much larger scale. Therefore, it is critical to design preventive and effective treatment strategies.

The control of malware spreading can be considered as an optimal control problem that defines a trade-off solution between the cost of fast and periodic development of patches and the value of the recovery of the devices. At the same time, information propagation of the vulnerability of the computing devices and personal accounts in social networks as well as the knowledge of the effective protection measures can help raise the awareness of the security threats and their solutions to reduce the number of infected devices. Generally, multiple types of viruses co-exist at the same time. Hence, we model the malware spreading as Susceptible-Infected-Recovered-Susceptible (SIRS) dynamics in which the population of devices is grouped into several subpopulations, i.e., the susceptible (S), the infected (I) and recovered (R). In addition, a group of infected nodes is also divided into several subgroups. The SIRS dynamics describe the evolution of the population size that can be controlled using special patching and recovery. Spreading of information is also described by the modified SIR model, which includes susceptible (S), warned (W), and recovered (R) nodes. Here, warned nodes are informed of the necessity of protection of their accounts and devices from their neighbors.

The goal of the work is to combine the two epidemic processes in one model. One epidemic process describes the dissemination of the information and the other one is the spreading of the viruses. We consider a generalized Susceptible-Warned-Infected-Recovered-Susceptible (SWIRS) model, which extends the model for information spreading by incorporating the SIRS model that describes the propagation of two types of malware. In the paper, we present the stability analysis of SWIRS model, formulate a controlled SWIRS model, and show the structure of the optimal policies of spreading information about virus protection and optimal treatment. Moreover, we carry out a series of numerical simulations to corroborate the results.

Recent literature has seen a surge of interest in using optimal control and stability equilibrium analysis to study malware protection in computer networks, social networks, and ad-hoc networks (See \cite{Fedyanin,Zuzek,Taynitskiy,Taynitskiy:2017,Junaid1,Junaid2,Moon,yunhan1, yunhan2}). Moreover, the clusters of the population play an important role. Several waves of the viruses propagation might occur due to sequential propagation information from one cluster to another even when a single cluster model might predict just a monotone spreading. 

In this paper, we establish a control-theoretic model to design optimal quarantining and immunization strategies to mitigate the impact of epidemics on our society. The recent spreading of ransomware (e.g., CryptoLocker, CryptoDefense, or CryptoWall) has spread using spam emails to extort money from home users and businesses alike by locking files on a PC or network storage (See \cite{Luo,wired}).
Mean-field dynamical systems are used to model the underlying evolution of the host subpopulations.  In \cite{Wang}, many variants of optimal control models of SIR-epidemics are investigated in the context of medical vaccination and health promotion campaigns. Previous studies have shown the application of epidemic frameworks to the models of network protection as in \cite{Mieghem,Sahneh,Vespignani,Junaid1,Junaid2,Taynitskiy:2017b,Taynitskiy:2018,Altman:2019}. Many different research works have provided many variants of epidemic models in computer security. Spreading information on social networks is presented in \cite{Moore}.

The rest of the paper is organized as follows. Section \ref{epidemic} presents the controlled SWIRS mathematical model. In Section \ref{Model formulation}, we formulate the SWIRS model. Sections \ref{Stability_analysis} and \ref{Lyapunov} show the stability analysis of the disease-free equilibrium. Section \ref{structure} describes the optimal control problem and Section \ref{Optimal_Control} presents the structure of optimal protection and information spreading policies. In Section \ref{Deterministic_network_cluster_model}, theoretical results are applied to the case of clusterized population. Section \ref{numerical} presents the series of numerical experiments. Section \ref{conclusion} concludes the paper.

\section{Deterministic population model}\label{epidemic}
\subsection{Model formulation}\label{Model formulation}
In this section, we formulate a two-level modified SIRS model (Susceptible-Infected-Recovered-Susceptible) with two different types of viruses circulated in a population of size $N$. This auxiliary partitioning allows capturing two processes that occur in both computer and social networks. The first process is the propagation of information on harmful malware attacks and the protection of personal data, documents, projects, etc. We consider this spreading process as the first level hierarchy in the Susceptible-Warned-Infected-Recovered-Susceptible (SWIRS) model. The second process, which corresponds to the physical propagation of antivirus software, is considered as the second level of the model, which is a modified Susceptible-Infected-Recovered-Susceptible (SIRS) model with two competitive viruses. Thereby, in contrast to classical SIRS models, where populations are divided into three groups: \textit{Susceptible} $(S)$, \textit{Infected }$(I)$, and \textit{Recovered} $(R)$, here the Infected subgroup is divided into two subgroups: a subgroup of nodes infected by the first type of virus $V_1$ and the subgroup infected by the second type $V_2$. Spreading information on the first level adds a new group \textit{Warned} $(W)$ into consideration. This group consists of the nodes, which have received information about possible the risks of virus attack/spreading and ready to use special tools for protections.

%We model the epidemic process as a system of nonlinear differential equations, where  $n_S$, $n_W$ $n_{V_1}$, $n_{V_2}$ and $n_R$ correspond to the number of susceptible, warned, infected and recovered nodes, respectively. The total number of nodes in the network during the entire process remains constant and equal to $n_S+n_W+n_{V_1}+n_{V_2}+n_R=N$. Let $S(t)=\frac{n_S(t)}{N}$, $W(t)=\frac{n_W(t)}{N}$, $I_1(t)=\frac{n_{V_1}(t)}{N}$, $I_2(t)=\frac{n_{V_2}(t)}{N}$, $R(t)=\frac{n_R(t)}{N}$  as  a fraction of the \textit{Susceptible}, the \textit{Warned}, the \textit{Infected}  and the \textit{Recovered} nodes, respectively.
\begin{figure}[h!]
\begin{center}
  \includegraphics[width=68mm]{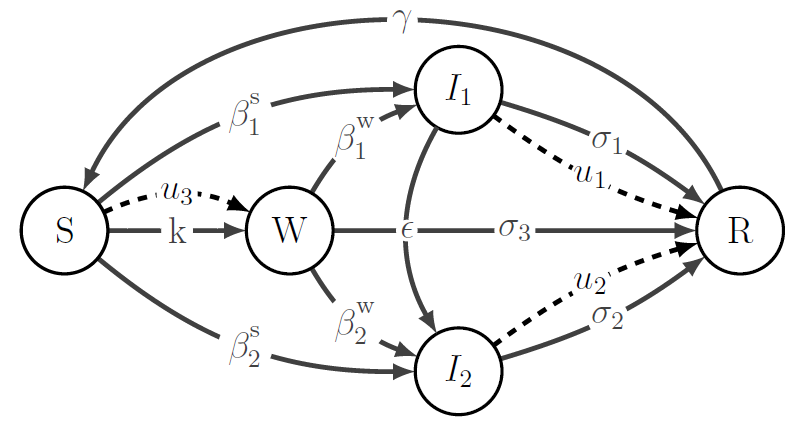}
  \caption{The scheme of transitions between groups $S$, $W$, $I_1$, $I_2$, $R$.}
\end{center}\label{fig:Epidemic_dynamics_2_clasters}
\end{figure}
We model the epidemic process as a system of nonlinear differential equations. The total number of nodes in the network during the entire process remains constant and equal to $n_S+n_W+n_{V_1}+n_{V_2}+n_R=N$. Let $S(t)=\frac{n_S(t)}{N}$, $W(t)=\frac{n_W(t)}{N}$, $I_1(t)=\frac{n_{V_1}(t)}{N}$, $I_2(t)=\frac{n_{V_2}(t)}{N}$, $R(t)=\frac{n_R(t)}{N}$  as  a fraction of the \textit{Susceptible}, the \textit{Warned}, the \textit{Infected},  and the \textit{Recovered} nodes, respectively. At the beginning of the epidemic, at time $t=0$, the majority of the individuals are in the Susceptible state, and a small fraction of individuals are infected by different types of virus. Hence, initial states are $S(0)=S^0>0,$ $W(0)=W^0\ge 0,$ $I_{1}(0)=I_{1}^0>0,$ $I_{2}(0)=I_{2}^0>0$ and  $R(0)=R^0=1-S^0-W^0-I_{1}^0-I_{2}^0.$

%\displaystyle{\frac{dW}{dt}}
Behaviour of the system is described by a system of nonlinear differential equations:
\begin{equation}\label{SIR1}\begin{array}{l}
dS/dt=-kWS-\beta^S_1 SI_1-\beta^S_2 SI_{2}+\gamma R-u_3S;\\
%\end{array}\end{equation}
%\begin{equation}\label{SIR1}\begin{array}{l}
dW/dt=kWS-\beta^W_1 WI_1-\beta^W_2 WI_{2}+u_3S-\sigma_3W;\\
dI_1/dt=\beta^S_1 SI_1+\beta^W_1WI_1-\varepsilon I_1I_2-\sigma_1I_1 -u_1I_1;\\
dI_2/dt=\beta^S_2 SI_2+\beta^W_2WI_2+\varepsilon I_1I_2-\sigma_2I_2 -u_2I_2;\\
dR/dt=\sigma_1I_1+u_1 I_1+\sigma_2I_2+u_2I_{2}+\sigma_3W-\gamma R,
\end{array}\end{equation}
where $\beta^S_i$ are infection rates for susceptible nodes for virus $V_i$, $i=1,2$ and $\beta^W_i$ are infection rates for the warned nodes. On the second level of the epidemic process, we can view a self-recovery rate $\sigma_1$ for virus $V_1$ or $\sigma_2$ for virus $V_2$ as the probability that infected nodes from subgroups $I_1$ or $I_2$ are recovered from the infection without incurring any costs on our system. On the first level, nodes that are informed of virus attacks have recovery rate $\sigma_3$. Without loss of generality, we can say that the second  virus  $V_2$ is stronger than the first $V_1$, and with the probability $\varepsilon$  virus $V_2$ can supersede the first virus in the node infected by the first virus.

The application of antivirus patches reduces the number of the infected. It can be interpreted as control parameters by $u_1(t)$ and $u_2(t)$ in (\ref{SIR1}), where $u_i$ are the fractions of the infected under treatment, $u_1(t),u_2(t)\in[0,1], \mbox{for all}\ t$. The warned nodes can avoid an epidemic by taking special quarantine measures. Control parameter $u_3(t)$ is the fraction of susceptible nodes that become warned of the virus spreading at time $t$.
\subsection{Stability analysis}\label{Stability_analysis}
In this section, the stability of the equilibrium points of the uncontrolled system is presented, where $u_i=0, i=1,2,3$ (\cite{Capasso,Allen,Wu,Sharma}). The disease-free equilibrium  is defined as the steady-state, where $I_1(t)=I_2(t)=0$ for any $t$. Assume that $I_1(t)=I_2(t)=0$, which means that the system is independent of the viruses, and we obtain simplified SWIRS-model:
{\small\begin{equation}\label{SIRdiseasefree}\begin{array}{l}
\dot{S}=-kWS+\gamma R;\\
\dot{W}=kWS-\sigma_3W;\\
\dot{I_1}=\dot{I_2}=0;\\
\dot{R}=\sigma_3W-\gamma R.\\
\end{array}\end{equation}}
By solving the system (\ref{SIRdiseasefree}), we obtain two disease-free equilibrium points:
{\small\begin{itemize}
    \item $E_1$: $S=1$, $W=I_1=I_2=R=0$;\\
    \item $E_2$: $\displaystyle S=\frac{\sigma_3}{k},\ W=\frac{\gamma(k-\sigma_3)}{k(\gamma+\sigma_3)}$, $I_1=I_2=0$, $\displaystyle R=\frac{\sigma_3(k-\sigma_3)}{k(\gamma+\sigma_3)}$.
\end{itemize}}
Local stability of the disease-free equilibrium points  is verified by studying   the real parts of eigenvalues $p_l$ of the Jacobian matrix at these points, i.e., $Re\ p_l\le 0$ for all $l$ (\cite{Capasso}). %Jacobian $Jac(S,W,I_1,I_2,R)$ for SWIRS-model is defined as follows:

1) Consider the first disease-free equilibrium point $E_1=(1,0,0,0,0)$. Define the Jacobian at this as $Jac_1$:
{\small
\begin{equation}
Jac_1=
\begin{pmatrix}
0 & -k & -\beta^S_1 & -\beta^S_2 & \gamma\\
0 & k-\sigma_3 & 0 & 0 & 0\\
0 & 0 & \beta^S_1-\sigma_1 & 0 & 0\\
0 & 0 & 0 & \beta^S_2-\sigma_2 & 0\\
0 & \sigma_3 & \sigma_1 & \sigma_2 & -\gamma
\end{pmatrix}.
\end{equation}}
This Jacobian has five eigenvalues $p_1=0$, $p_2=k-\sigma_3$, $p_3=\beta^S_1-\sigma_1$, $p_4=\beta^S_2-\sigma_2$, $p_5=-\gamma$.

\begin{prop1}
Since all parameters ($\beta^S_1,\beta^S_2,\sigma_1,\gamma,etc.$) are non-negative, equilibrium point $E_1$ will be asymptotically stable if the following conditions are hold:
\begin{equation}
k\le\sigma_3,\hskip 15pt \beta^S_1\le\sigma_1,\hskip 15pt \beta^S_2\le\sigma_2.
\end{equation}
\end{prop1}
\begin{figure}[h!]
\begin{center}
 \includegraphics[width=71mm]{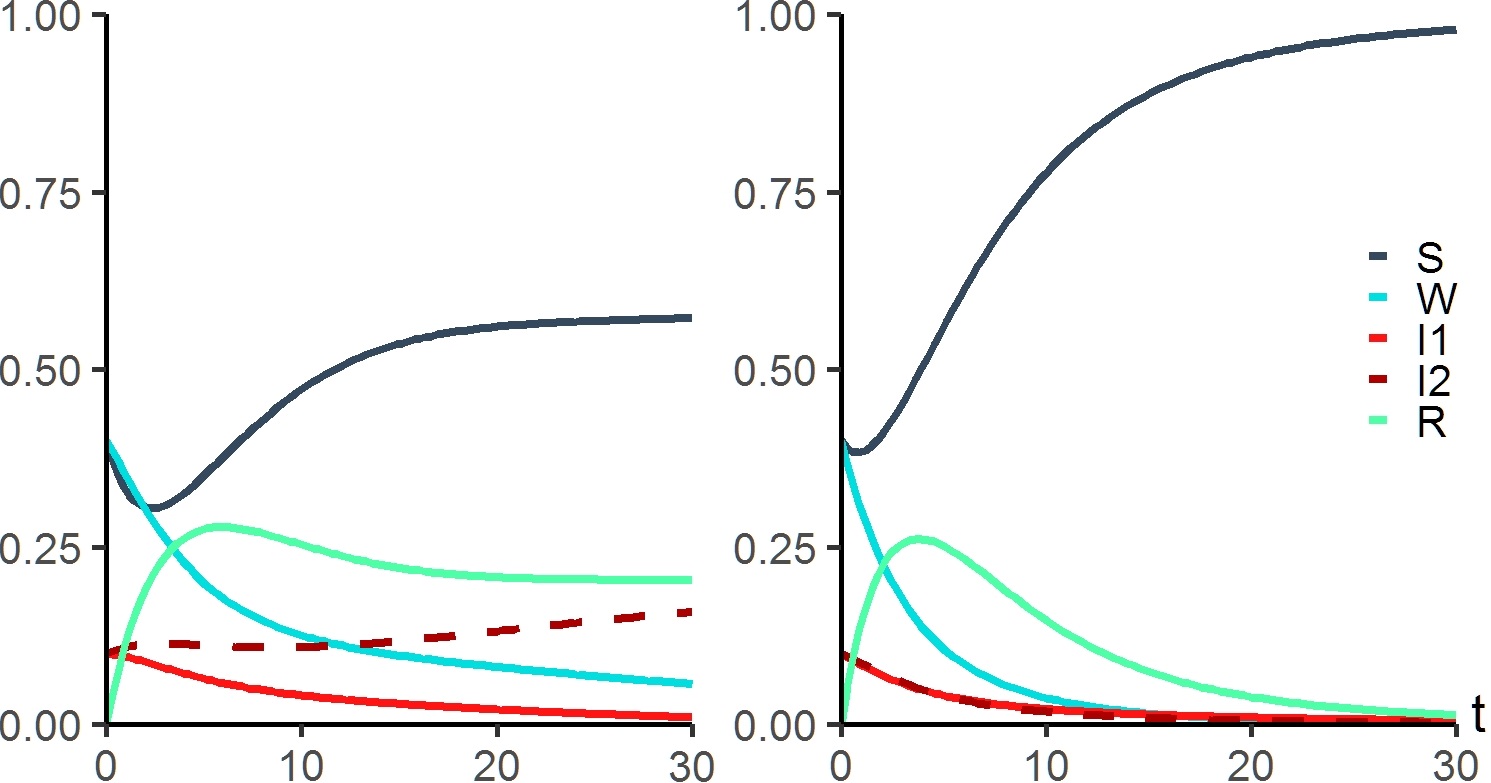}
\caption{Disease-free equilibrium point $E_1$: conditions (4) are not satisfied (left) and are satisfied (right).}
\end{center}\label{fig:DFE}
\end{figure}
2) For the second disease-free equilibrium point $E_2$, the Jacobian $Jac_2$ has the following form 
{\small\begin{equation}\label{jac2}
 \begin{pmatrix}
 kx_1 & -\sigma_3 & -x_2 & -x_3 & \gamma\\

 -kx_1 & 0 & \beta^W_1 x_1 & \beta^W_2 x_1 & 0\\

 0 & 0 & x_2-\beta^W_1x_1-\sigma_1 & 0 & 0\\

 0 & 0 & 0 & x_3-\beta^W_2x_1-\sigma_2 & 0\\

 0 & \sigma_3 & \sigma_1 & \sigma_2 & -\gamma
 \end{pmatrix},
 \end{equation}}
where $x_1=\frac{\gamma(\sigma_3-k)}{k(\gamma+\sigma_3)}$, $x_2=\frac{\beta^S_1\sigma_3}{k}$ and $x_3=\frac{\beta^S_2\sigma_3}{k}$.
% or

% \begin{equation}\label{jac2}
% Jac_2=
% \begin{pmatrix}
% \frac{\gamma(\sigma_3-k)}{\gamma+\sigma_3} & -\sigma_3 & -\frac{\beta^S_1\sigma_3}{k} & -\frac{\beta^S_2\sigma_3}{k} & \gamma\\

% -\frac{\gamma(\sigma_3-k)}{\gamma+\sigma_3} & 0 & \frac{\beta^W_1\gamma(\sigma_3-k)}{k(\sigma_3+\gamma)} & \frac{\beta^W_2\gamma(\sigma_3-k)}{k(\sigma_3+\gamma)} & 0\\

% 0 & 0 & \frac{\beta^S_1\sigma_3}{k}-\frac{\beta^W_1\gamma(\sigma_3-k)}{k(\sigma_3+\gamma)}-\sigma_1 & 0 & 0\\

% 0 & 0 & 0 & \frac{\beta^S_2\sigma_3}{k}-\frac{\beta^W_2\gamma(\sigma_3-k)}{k(\sigma_3+\gamma)}-\sigma_2 & 0\\

% 0 & \sigma_3 & \sigma_1 & \sigma_2 & -\gamma
% \end{pmatrix}.
% \end{equation}

Jacobian matrix $Jac_2$ has five eigenvalues:
\begin{itemize}
    \item $\displaystyle p_{q}=-\beta^W_qx_1-\sigma_q+\beta^S_q\sigma_3/k$;
    \item $\displaystyle p_{3,4}=\frac{\gamma(\gamma+k-2\sigma_3)\pm \sqrt D}{2(\sigma_3+\gamma)}$;
    \item $p_5=0$,
\end{itemize}

where $q\in \{1,2\}$ and $$D=\gamma^2(\gamma+k-2\sigma_3)^2+4\gamma(k-\sigma_3)(\gamma+\sigma_3)^2.$$
The following conditions define  the stability disease-free equilibrium.

\begin{prop1}
Equilibrium point $E_2$ will be asymptotically stable if the following conditions hold 

\begin{itemize}

    \item If  $D\le0$, then conditions are
$$\gamma\sigma_3(\beta^W_q-\beta^S_q)+k\gamma(\sigma_q-\beta^W_q)+\sigma_3(k\sigma_q-\sigma_3\beta^S_q)\ge0,$$
$$\gamma+k\le 2\sigma_3,\ \ \  q\in\{1,2\}.$$
    \item If $D>0$, then conditions are
$$\gamma\sigma_3(\beta^W_q-\beta^S_q)+k\gamma(\sigma_q-\beta^W_q)+\sigma_3(k\sigma_q-\sigma_3\beta^S_q)\ge0,$$
$$\gamma(\gamma+k-2\sigma_3)\pm \sqrt D\le0,\ \ \  q\in\{1,2\}.$$
\end{itemize}
\end{prop1}
\subsection{Global stability analysis of disease-free equilibrium $E_0$}\label{Lyapunov}
In this subsection, the global stability of disease-free equilibrium  $E_0(S,0,0,0,R) $ is discussed. For this purpose, we use the following Lyapunov function:
\begin{equation}
L(W, I_1, I_2)=W+I_1+I_2.
\end{equation}
here function $L(0)=0$ and $L(\cdot)\ge 0$ otherwise. The derivative of $L(W, I_1, I_2)$ with respect to the system (\ref{SIR1}) gives:
\begin{equation}\label{Lyapunov_func}\begin{array}{l}
L(W, I_1, I_2)|_{(\ref{SIR1})}=\dot{W}+\dot{I_1}+\dot{I_2}=(kS-\sigma_3) W+\\\hskip 70 pt (\beta_1^S S-\sigma_1) I_1+(\beta_2^S S-\sigma_3) I_2.
\end{array}\end{equation}
The disease-free equilibrium point is asymptotically stable if  the derivative $L(\cdot)|_{(\ref{SIR1})}<0$. This condition is satisfied if  the following conditions hold:
\begin{equation}\label{lyapunov_func_cond}
(kS-\sigma_3)<0,\ (\beta_1^S S-\sigma_1)<0,\ (\beta_2^S S-\sigma_3)<0,
\end{equation}
since variables $W$, $I_1$, $I_2$ are nonnegative. Conditions (\ref{lyapunov_func_cond}) show that if the self recovery rates are higher than the infection rates, then the epidemics vanishes. %Furthermore $\frac{d L}{dt}
%= 0$ if and only if $I_i = 0, i=1,2$, $W=0$. Thus, the largest compact invariant
%set in $\{(S,W, I_1, I_2,R) РѓС‘ С“РЋ : .L = 0}, when R0 . 1, is the singleton {E0}.
\subsection{Optimal Control of Epidemics}\label{structure}
%In Section I, we show the necessity of application of controls to protect the network from the virus attacks. 
It is clear that the protection measures have their costs. Let the objective function $J$ be the sum of two functionals, which correspond to the two levels of the model. On the first level, functional $J_1$ describes the costs of the quarantine measures, i.e. the costs of disseminating information about the epidemics to susceptible nodes. On the second level, functional $J_2$ defines the cost of antivirus treatment and includes the costs incurred by infected nodes, costs of spreading antivirus, and the benefit from the recovered nodes.

%Consider $S(t)=\frac{\sum_{j=1}^{N}S_j(t)}{N}$, $W(t)=\frac{\sum_{j=1}^{N}W_j(t)}{N}$, $I_1(t)=\frac{\sum_{j=1}^{N}{I_{1j}}(t)}{N}$, $I_2(t)=\frac{\sum_{j=1}^{N}{I_{2j}}(t)}{N}$, $R(t)=\frac{\sum_{j=1}^{N}R_j(t)}{N}$  as  a fraction of the \textit{Susceptible}, the \textit{Warned}, the \textit{Infected}  and the \textit{Recovered} nodes, respectively. Condition $S(t)+W(t)+I_1(t)+I_2(t)+R(t)=1$ holds for any $t$.
%

%Application of antivirus patches reduces the number of the infected individuals in the computer network and it is interpreted as control parameters which are  defined as $u_1$ and $u_2$ in (\ref{SIR1}). Where $u_i$ are fractions of the infected which are under treatment, $0\leq u_1(t)\leq 1$, $0\leq u_2(t)\leq 1, \mbox{for all}\ t$.  Warned nodes can avoid an epidemic by taking special quarantine measures. Control parameter $u_3(t)$  is responsible for the fraction of susceptible nodes that become warned about the virus spreading by special agents at time moment $t$.

At any given $t$,  $f_1(I_1(t)), f_2(I_2(t))$  are infection costs; $L(W(t))$ is the utility of the warned nodes. Function  $g(R(t))$  defines the benefit rate for recovered nodes; functions $h_1(u_1(t)), h_2(u_2(t))$ are costs for antivirus treatments and $h_3(u_3(t))$ is cost of information spreading. Here functions $f_i(I_i)$ are non-decreasing and twice-differentiable, convex functions, $f_i(0)=0$, $f_i(I_i)>0$ for $I_i >0,\ i=1,2$, $g(R)$ and $L(W)$ are non-decreasing and differentiable functions, and $h_i(u_i(t))$ is twice-differentiable and increasing function in $u_i(t)$ such as $h_i(0)=0$, $h_i(x)>0,\ \ i=1,2,3,$ when $u_i>0$. Also costs of information spreading are lower than costs for antivirus treatments $h_3(\cdot)<h_1(\cdot)$ and $h_3(\cdot)< h_2(\cdot)$.

The aggregated system costs over the time interval $[0,T]$ are defined as $J=J_1+J_2$,
%\begin{equation}\label{functional_J}\begin{array}{c}
%J=J_1+J_2, \end{array}
%\end{equation}
where
\begin{equation}\label{functional_SIR}\begin{array}{c}
J_1=\int_0^T h_3(u_3(t))-L(W(t)) dt,\\
J_2=\int_0^T
\sum^2_{q=1}\Big(f_q(I_q(t))+h_q(I_q(t))\Big)-g(R(t)).\end{array}
\end{equation}
and the optimal control problem is to minimize these costs, i.e.,
$\min_{\{u_1, u_2, u_3\}} J.$

By  using  Pontryagin's maximum principle (\cite{Pontryagin}), we  construct the optimal control $u(t)=(u_1(t), u_2(t), u_3(t))$ to the problem described above in Section \ref{epidemic}. To simplify the presentation, we use short-hand notations $S,I_1,u_1,$ etc. in place of $S(t),I_1(t),u_1(t),$ etc. Define the associated Hamiltonian $H$ and adjoint functions $\lambda_{S}(t)$, $\lambda_{W}(t)$, $\lambda_{I_1}(t)$, $\lambda_{I_2}(t)$, $\lambda_{R}(t)$ as follows:
%Hamiltonian (\textit{with mutation equation}):
\begin{equation}\label{Hamiltonian}
\begin{array}{l}
H=-f_1(I_1)-f_2(I_2)+g(R)-h_1(u_1)-h_2(u_2)+L(W)-\\\hskip 20pt 
h_3(u_3)+
(\lambda_W-\lambda_S)kWS+(\lambda_R-\lambda_W)\sigma_3W+ \\
\hskip 20pt(\lambda_{I_1}-\lambda_{S})\beta^S_1SI_1+
(\lambda_{I_2}-\lambda_{S})\beta^S_2SI_2+\\
\hskip 20pt (\lambda_{I_1}-\lambda_{W})\beta^W_1WI_1+
(\lambda_{I_2}-\lambda_{W})\beta^W_2WI_2+ \\
\hskip 20pt(\lambda_{I_2}-\lambda_{I_1})\varepsilon I_1I_2+(\lambda_{S}-\lambda_{R})\gamma R+
\\\hskip 20pt (\lambda_{R}-\lambda_{I_1})(\sigma_1+u_1)I_1+(\lambda_{R}-\lambda_{I_2})(\sigma_2+u_2)I_2+\\
\hskip 20pt (\lambda_W-\lambda_{S})u_3S.
\end{array}
\end{equation}
The adjoint system is defined as follows:
\begin{equation}\label{adjoint system}
\begin{array}{l}
\dot{\lambda}_S(t) = (\lambda_S-\lambda_W)kW+(\lambda_S-\lambda_{I_1})\beta^S_1I_1+\\
\hskip 37pt (\lambda_S-\lambda_{I_2})\beta^S_2I_2+(\lambda_S-\lambda_{W})u_3;\\
\dot{\lambda}_{W}(t) =-L'(W)+(\lambda_S-\lambda_W)kS+(\lambda_W-\lambda_{I_1})\beta^W_1I_1+\\
\hskip 37pt (\lambda_W-\lambda_{I_2})\beta^W_2I_2+(\lambda_W-\lambda_R)\sigma_3;\\
\dot{\lambda}_{I_1}(t)= f_1'(I_1)+(\lambda_S-\lambda_{I_1})\beta^S_1S+(\lambda_W-\lambda_{I_1})\beta^W_1W+\\
\hskip 37pt (\lambda_{I_1}-\lambda_{I_2})\varepsilon I_2+(\lambda_{I_1}-\lambda_{R})(\sigma_1+u_1);\\
%\end{array}
%\end{equation}
%
%\begin{displaymath}\begin{array}{l}
\dot{\lambda}_{I_2}(t)= f_2'(I_2)+(\lambda_S-\lambda_{I_2})\beta^S_2S+(\lambda_W-\lambda_{I_2})\beta^W_2W+\\
\hskip 37pt (\lambda_{I_1}-\lambda_{I_2})\varepsilon I_1+(\lambda_{I_2}-\lambda_{R})(\sigma_2+u_2); \\

\dot{\lambda}_R(t)=-g'(R)+(\lambda_R-\lambda_{S})\gamma,
\end{array}%\end{displaymath}
\end{equation}
with the transversality conditions given by
\begin{equation}\label{transversality}
\lambda_{S}(T)=\lambda_{W}(T)=\lambda_{I_1}(T)=\lambda_{I_2}(T)=\lambda_{R}(T)=0.
\end{equation}

According to Pontryagin's  maximum principle, there exist continuous and piece-wise
continuously differentiable co-state functions $\lambda_r(t),\ r\in\{S,W,I_1,I_2,R\}$  that satisfy (\ref{adjoint system}) and (\ref{transversality}) for  $t\in [0, T]$ together with continuous functions $u^*_1(t)$, $u^*_2(t)$ and $u^*_3(t)$:
\begin{equation}\label{maximum_principle}\begin{array}{l}
(u^*_1, u^*_2, u^*_3)\in \\\hskip 15pt \textrm{arg} \max\limits_{{u}_1, {u}_2,{u}_3\in
[0,1]} H(\lambda,S,W, I_1, I_2, R,{u}_1,
{u}_2,{u}_3).\end{array}
\end{equation}
\subsection{Structure of Optimal Control}\label{Optimal_Control}
In this subsection, we construct the structure of the optimal control $u^*(t)=(u^*_1(t), u^*_2(t),u^*_3(t))$.
\begin{prop1}\label{proposition_1}
The following statements hold for the optimal control problem described in Section \ref{epidemic}:

\begin{itemize}
  \item When $h_i(\cdot)$ are concave functions, then there exists  $t_0\in [0,T]$ such that for any $i=1,2,3:$

  $$u^*_i(t)=\left\{\begin{array}{l} 1,\ \mbox{for}\ 0\le t\le t_0;\\
 0,\ \mbox{for}\ t_0<t\le T.\end{array}\right.$$
  \item When $h_i(\cdot)$ are strictly convex functions, then there exist the time $t_0, t_1$, $0<t_0<t_1<T$ such that for any $i=1,2,3$ ($\alpha(t)\in(0,1)$):
$$
u^*_i(t)= \left\{\begin{array}{l}1, \hskip 23pt 0\le t\le t_0;\\
 \alpha(t),\ \ \ t_0<t\le t_1;\\
 0,\hskip 23pt t_1<t\le T.\\
\end{array}\right.$$
\end{itemize}
\end{prop1}

We define functions $\varphi_i(t)$ as follows:
\begin{equation}\label{Hamiltonian1}
\begin{array}{l}
\varphi_q(t)=(\lambda_{R}(t)-\lambda_{I_q}(t))I_{q}(t),\ q\in \{1,2\}, \\ 
\varphi_3(t)=(\lambda_{W}(t)-\lambda_{S}(t))S(t).
\end{array}
\end{equation}

 To prove Proposition 1, we consider the following auxiliary lemma.

\begin{prop1}\label{prop_0} Functions $\varphi_i, \ i=\overline{1,3}$ are decreasing functions of\ $t$ for  $t\in [0, T].$
\end{prop1}

Let's rewrite the Hamiltonian in terms of function $\varphi_i(t)$:
{\small\begin{equation}\label{Hamiltonian1_1}
\begin{array}{l}
H=-f_1(I_1)-f_2(I_2)+g(R)+L(W)+(\lambda_W-\lambda_S)kWS+\\
\hskip 20pt (\lambda_{I_1}-\lambda_{S})\beta^S_1SI_1+
(\lambda_{I_2}-\lambda_{S})\beta^S_2SI_2+\\
\hskip 20pt (\lambda_{I_1}-\lambda_{W})\beta^W_1WI_1+
(\lambda_{I_2}-\lambda_{W})\beta^W_2WI_2+\\
\hskip 20pt (\lambda_{I_2}-\lambda_{I_1})\varepsilon I_1I_2+
(\lambda_{R}-\lambda_{I_1})\sigma_1I_1+ (\lambda_{R}-\lambda_{I_2})\sigma_2I_2+\\
\hskip 20pt (\lambda_{S}-\lambda_{R})\gamma R+
(-h_1(u_1)+\varphi_1u_1)+\\
\hskip 20pt (-h_2(u_2)+\varphi_2u_2)+(-h_3(u_3)+\varphi_3u_3).
\end{array}
\end{equation}}
We can divide this maximization problem into three subproblems and find optimal control $u^*_1(t), u^*_2(t)$ and $u^*_3(t)$, separately.
\begin{equation}\label{maximum_principle_3}\begin{array}{l}
\max\limits_{u_1, u_2,u_3}
[-h_1(u_1)+\varphi_1 u_1-h_2(u_2)+\varphi_2
u_2-\\
\hskip 18pt h_3(u_3)+\varphi_3u_3]=\max\limits_{u_1}[-h_1(u_1)+\varphi_1
u_1]+\\
\hskip 17pt \max\limits_{u_2}[-h_2(u_2)+\varphi_2 u_2]+\max\limits_{u_3}[-h_3(u_3)+\varphi_3 u_3].
\end{array}\end{equation}
We obtain the following derivatives:
\begin{equation}\label{dH_du}\begin{array}{l}
\displaystyle\frac{\partial H}{\partial u_i}=-\dot{h}_i(u_i)+\psi_i=0, \ i=\overline{1,3}.
\end{array}\end{equation}
As $h_i(u_i)$ are increasing functions and $I_q\ge0$ and $S\ge0$, then the Hamiltonian reaches its maximum if $
\psi_i=\dot{h}_i(u_i)\ge0, \ i=1, 2 , 3.$ We can find such $u_i$ if and only if the following  conditions are satisfied:
 $\lambda_R(t)-\lambda_{I_1}(t)\ge0$, $\lambda_R(t)-\lambda_{I_2}(t)\ge0$ and $\lambda_W(t)-\lambda_{S}(t)\ge0$. To complete the proof of proposition, we consider the auxiliary lemma.

\begin{lemma1}\label{lemma_1}
 For all $t\in[0,T]$, we have $\lambda_R(t)-\lambda_{I_1}(t)\ge0$, $\lambda_R(t)-\lambda_{I_2}(t)\ge 0$ and $\lambda_W(t)-\lambda_{S}(t)\ge0$.
\end{lemma1}

Proof of the Lemma \ref{lemma_1} is based on the following  properties:
{{\bf Property 1:} Let $v(t)$ be a continuous and piece-wise differential function of $t$. Let $v(t_1) = L$ and $v(t) > L$ for all $t\in(t_1, \ldots , t_0]$. Then $\dot{v}(t_1^+)\ge0$. Where $v(t_1^+)=\lim\limits_{x\to t_1+0} v(x)$.} \\
{{\bf Property 2:} For any convex and differentiable function $y(x)$, which is $0$ at $x = 0$, $y'(x)x-y(x)\ge 0$ for all $x\ge 0$.}

%\begin{proper}\label{property_1}
%Let $v(t)$ be a continuous and piece-wise differential function of $t$. Let $v(t_1) = L$ and $v(t) > L$ for all $t\in(t_1, \ldots , t_0]$. Then $\dot{v}(t_1^+)\ge0$. Where $v(t_1^+)=\lim\limits_{x\to 0} v(x)$.
%\end{proper}

%\begin{proper}\label{property_2}
%For any convex and differentiable function $y(x)$, which is $0$ at $x = 0$, $y'(x)x-y(x)\ge 0$ for all $x\ge 0$.
%\end{proper}

We divide our proof into two parts. In the first part, we consider the case when $t=T$ and show that derivatives of the functions $\lambda_R(t)-\lambda_{I_1}(t)$, $\lambda_R(t)-\lambda_{I_2}(t)$ and $\lambda_W(t)-\lambda_S(t)$ are less or equal to zero to show that they are non-increasing. In the second part, we use the method of proof by contradiction and show that on the whole interval $[0,T]$ these functions are also non-negative.

\textbf{Step I.} At time  $T$, according to (\ref{transversality}), we have that $\ \lambda_R(T)-\lambda_{I_1}(T)=0$, $\ \lambda_R(T)-\lambda_{I_2}(T)=0$ and
$\\ \lambda_W(T)-\lambda_S(T)=0.$ From (\ref{adjoint system}) it is obtained that all he derivatives are non-positive
\begin{equation}\begin{array}{l}
\dot{\lambda}_R(T)-\dot{\lambda}_{I_q}(T)=-\dot{g}(R(T))-\dot{f}_q(I_q(T)\le0,\ q\in \{1,2\}, \\
\dot{\lambda}_W(T)-\dot{\lambda}_S(T)=-\dot{L}(W(T))\le0.
\end{array}\end{equation}
Since $g(\cdot)$, $f_1(\cdot)$, $f_2(\cdot)$ and $L(\cdot)$ are increasing functions, at  time  $T$ all functions are equal to $0$ and their derivatives are less or equal to $0$, then we can obtain that $\lambda_R(t)-\lambda_{I_1}(t)$, $\lambda_R(t)-\lambda_{I_2}(t)$ and $\lambda_W(t)-\lambda_S(t)$ are non-increasing functions at $t=T$.

\textbf{Step II.} In this step, we   show by contradiction that $\lambda_R(t)-\lambda_{I_1}(t)\ge0$ for all $t\in[0,T]$. Proofs for the $\lambda_R(t)-\lambda_{I_2}(t)$ and $\lambda_W(t)-\lambda_S(t)$ use the same method and we will leave it to the readers.

The system of ODE (1) is autonomous, and, hence, the Hamiltonian  and the control do not depend on the variable independent $t$. From (\ref{Hamiltonian}), we obtain
\begin{equation}\label{ham_eq}\begin{array}{l}
H+f_2+\sum _i h_i-g(R)-L(W)-(\lambda_{I_2}-\lambda_{S})\beta^S_2SI_2-\\
(\lambda_R-\lambda_W)\sigma_3W-(\lambda_{I_2}-\lambda_{W})\beta^W_2WI_2-\\
(\lambda_{R}-\lambda_{I_2})(\sigma_2+u_2)I_2- (\lambda_{S}-\lambda_{R})\gamma R-\\
(\lambda_W-\lambda_{S})(kW+u_3)S\le -f_1(I_1(T))\le 0.
\end{array}\end{equation}
Suppose that there exists time moment $t^*\in(0,T)$ at which $\lambda_R(t^*)-\lambda_{I_1}(t^*)=0$. Using (\ref{ham_eq}), consider the derivative of this function at the time moment $t^{*+}$:
{\begin{equation}\label{lemma_eq}\begin{array}{l}
\dot{\lambda}_{R}(t^{*+}) -\dot{\lambda}_{I_1}(t^{*+})=-\dot{f_1}(I_1)-\dot{g}(R)-(\lambda_S-\lambda_R)\gamma+\\
(\lambda_{I_1}-\lambda_S)\beta^S_1S+(\lambda_{I_1}-\lambda_W)\beta^W_1W+(\lambda_{I_2}-\lambda_{I_1})\varepsilon I_2+\\
(\lambda_R-\lambda_{I_1})(\sigma_1+u_1)I_1.
\end{array}\end{equation}}
%\frac{1}{I_1}(H+f_1+f_2+h_1+h_2+h_3-g(R)-L(W)-\\
%(\lambda_W-\lambda_S)kWS-(\lambda_R-\lambda_W)\sigma_3W-(\lambda_{I_1}-\lambda_{S})\beta^S_1SI_1-\\
%(\lambda_{I_2}-\lambda_{S})\beta^S_2SI_2-(\lambda_{I_1}-\lambda_{W})\beta^W_1WI_1-(\lambda_{I_2}-\lambda_{W})\beta^W_2WI_2-\\
%(\lambda_{I_2}-\lambda_{I_1})\varepsilon I_1I_2-(\lambda_{R}-\lambda_{I_2})\sigma_2I_2- (\lambda_{S}-\lambda_{R})\gamma R-\\
%(\lambda_R-\lambda_{I_2})u_2I_2- (\lambda_W-\lambda_{S})u_3S)
From (\ref{Hamiltonian}) and (\ref{lemma_eq}), we can obtain that
{\begin{equation}\label{lemma_eq}\begin{array}{l}
\dot{\lambda}_{R}(t^{*+}) -\dot{\lambda}_{I_1}(t^{*+})=-\frac{1}{I_1}(\dot{f_1}I_1-f_1)-(\lambda_S-\lambda_R)\gamma-\\
\dot{g}(R)+\frac{1}{I_1}\big(H+f_2+h_1+h_2+h_3-g(R)-L(W)-\\
(\lambda_R-\lambda_W)\sigma_3W-(\lambda_{I_2}-\lambda_{S})\beta^S_2SI_2-\\
(\lambda_{I_2}-\lambda_{W})\beta^W_2WI_2-(\lambda_{R}-\lambda_{I_2})(\sigma_2+u_2)I_2- \\
(\lambda_{S}-\lambda_{R})\gamma R-(\lambda_W-\lambda_{S})(kW+u_3)S\big).
\end{array}\end{equation}}
Here, $f_1(I_1)$ is convex and differentiable function, from Property 2 and  (\ref{ham_eq}) we obtained that $\dot{\lambda}_{R}(t^{*+})-\dot{\lambda}_{I_1}(t^{*+})\le 0$, but according to our assumption the derivative should be greater or equal to zero. That leads to contradiction and completes the proof that $\lambda_R(t)-\lambda_{I_1}(t)\ge0$ for all $t\in[0,T]$. Using the same method, we can prove that $\lambda_R(t)-\lambda_{I_2}(t)\ge0$ and $\lambda_W(t)-\lambda_S(t)\ge0$ for $t\in[0,T].$

Functions $\lambda_R(t)-\lambda_{I_1}(t)$, $\lambda_R(t)-\lambda_{I_2}(t)$ and $\lambda_W(t)-\lambda_S(t)$ are non-negative at the interval $[0,T]$ and at  $t=T$ the derivatives of these functions are less or equal to 0 that completes the proof of the Proposition \ref{prop_0}.

%Proof of Proposition 1 is similar to the proof in (\cite{Taynitskiy}). From this proof, we obtain conditions for the optimal control $u^*_1(t), u^*_2(t)$ and $u^*_3(t)$ that depend on the properties of functions $h_i(u_i(t))$, $i=1,2,3$.

\subsubsection{ Functions $h_i(\cdot)$ are concave}$\\$
Let $h_i(\cdot)$ be a concave functions ($h''_i(\cdot)<0$), then according to (\ref{Hamiltonian}) Hamiltonian is a convex function of $u_i,\ i=\overline{1,3}$. There could be two different options for $u_i\in[0,1]$ that maximimize the Hamiltonian. If $-h_i(0)+\varphi_i\cdot 0>-h_i(1)+\varphi_i\cdot 1$ or $h_i(1)>\varphi_i$, then optimal control -- $u_i=0$ (Fig. 3 (right)), otherwise -- $u_i=1$ (Fig. 3 (left)). For $i=\overline{1,3}$, the optimal control parameters $u_i(t)$ are defined as follows:
%There can be at most one time moment $t\in[0,T]$ at which $\varphi_i(t)=h'_i(u_i(t))$, moreover if such moment exists, for example, $t_1$, then $\varphi_i(t)<h_i(1)$ on $0\leq t< t_1$ and $\varphi_i(t)\geq h_i(1)$ on $t_1\leq t< T$. There could be only two unique points that maximises the Hamiltonian and it is either in $0$ or $1$.  
\begin{figure}[h!]
\begin{center}
  \includegraphics[width=56mm]{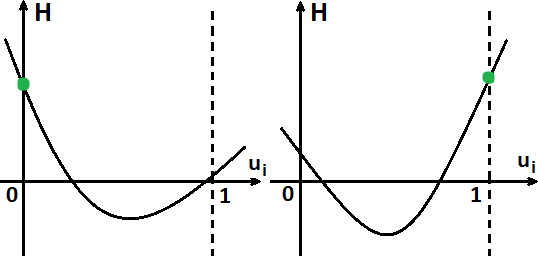}
  \caption{Hamiltonian when functions $h_i(\cdot)$ are concave.}
  \end{center}
\end{figure}
\begin{equation}\label{control}
u^*_i(t)=\left\{\begin{array}{l} 0,\ \varphi_i(t)<h_i(1),\\
1, \ \varphi_i(t)\geq h_i(1).\\
\end{array}\right.
\end{equation}

\subsubsection{Functions $h_i(\cdot)$ are strictly convex}$\\$
Let $h_i(\cdot)$ be a strictly convex functions ($h''_i(\cdot)>0$), then Hamiltonian is concave function. Consider the following derivative:
{\small\begin{equation}\label{der}\begin{array}{l}
\frac{\partial}{\partial x}(-h_i(x)+\varphi_i x)\mid_{x=x_i}=0,
\end{array}\end{equation}}
where $x\in [0,1]$, $u^*_i(t)=x_i$. There could be three different types of points at which the Hamiltonian reaches its maximum (Fig. 4). To find them, we need to consider the derivatives of the Hamiltonian at $u_i=0$ and $u_i=1$. If the derivatives (\ref{der}) at $u_i=0$ are non-increasing ($-h'_i(0)+\varphi_i\le0$), then the value of the control that maximizes the Hamiltonian is less than 0, and according to our restrictions ($u_i\in [0,1]$) optimal control will be equal to 0 (Fig. 4a). If the derivatives at $u_i=1$ are increasing ($-h'_i(1)+\varphi_i>0$), it means that the value of the control that maximizes the Hamiltonian is greater than 1. Hence the optimal control will be equal to 1 (Fig. 4c), otherwise, we can find such value $u^*_i\in (0,1)$ (Fig. 4b):
\begin{figure}%[h!]
\begin{center}
  \includegraphics[width=83mm]{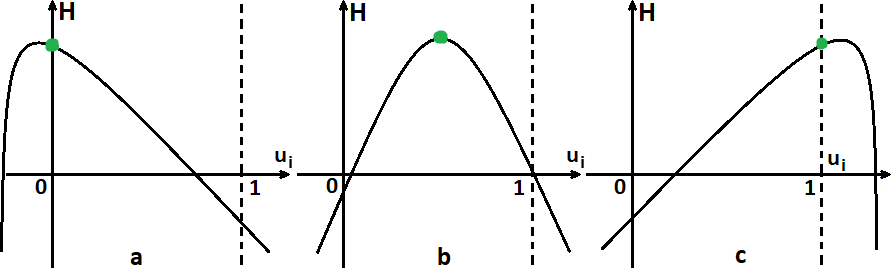}
  \caption{Hamiltonian when functions $h_i(\cdot)$ are convex.}
\end{center}
\end{figure}
\begin{equation}\label{control1}
u^*_i(t)=\left\{\begin{array}{l} 0, \ \ \varphi_i\leq h'_i(0),\ i=1,2,3;\\
h'^{-1}(\varphi_i), \ \ \ h'_i(0)<\varphi_i\leq h'_i(1),\ i=1,2,3.\\
1,\ \ \ h'_i(1)<\varphi_i,\ i=1,2,3.
\end{array}\right.
\end{equation}
Functions $\varphi_i(t)$, $h_i'(t)$, $u^*_i(t)$ are continuous at all $t\in [0, T].$ In this case $h_i$ is strictly convex and $h_i'$ is strictly increasing functions, so $h'(0)<h'(1)$. Thus there exist  points $t_0$ and $t_1$ $(0<t_0<t_1<T)$ so that conditions (\ref{control1}) are satisfied, and according to $\varphi_i$ are decreasing functions.

After obtaining the optimal control $u^*_1(t)$ and $u^*_2(t)$, we need to sort all infected nodes by the number of neighbors and treat them in order, starting with the first one on the list. Similar procedure is used to find the number of susceptible nodes, among which it is necessary to disseminate information about virus attacks, using the structure of the $u^*_3(t)$.

\section{SWIRS model on meta-population network }\label{Deterministic_network_cluster_model}
The clustering of the nodes in the network can be considered as a natural extension of the SWIRS model from Section \ref{epidemic}.  We assume that all nodes inside the one cluster follow the same behavioral rules. However, the infection can be transferred among  clusters. For this reason, we consider a case of a network with $N$ nodes, which can be divided into several clusters. Here, the matrix $A=\{a_{\tau \mu}\}$ is the adjacency matrix of the first level of SWIRS model, where information about possible consequences of malware attacks is spreading, and $B=\{b_{\tau \mu}\}$ is the adjacency matrix of the second level, where special antivirus patches are applied. Denote as $ka_{\tau \mu}$ the probability that a node from cluster $\tau$ of size $N_\tau$ and a node from a cluster $\mu$ of size $N_\mu$ change their states from $S$ to $W$ at every time instant. The probability that a susceptible node from cluster $\tau$ will be infected due to the contact with a node from a cluster $\mu$, infected by virus $V_l, l \in \{1,2\}$ is equal to $\beta^S_{V_l}b_{\tau \mu}$. A warned node from cluster $\tau$ will be infected by virus $V_l,\ l=\{1,2\}$ through the contact with the node from a cluster $\mu$ with  probability  $\beta^W_{V_l}b_{\tau \mu}$.

Vector $X_j(t)=(S_j(t),W_j(t),I_{1j}(t),I_{2j}(t),R_j(t))$ defines the proportions distribution of being in each of the states for the cluster $j=1,\ldots, M$ at $t$. For any $t\in [0,T]$, the sum of the probabilities for any node $j$ is equal to  $S_j(t)+W_j(t)+{I_{1j}}(t)+{I_{2j}}(t)+R_j(t)=1$. All other parameters in the system remain the same as in Section 3.1.
%\displaystyle{\frac{dS_j(t)}{dt}}
This simultaneous process of information spreading and  patching is described by a system of nonlinear differential equations:
{\small\begin{equation}\label{SIR2}
 \begin{array}{l}
 dS_j(t)/dt= -kS_j(t) \sum_l{a_{jl}W_l(t)}-\beta^S_1 S_j(t) \sum_l{b_{jl}I_{1l}(t)}-\\
\hskip 38pt\beta^S_2 S_j(t) \sum_l{b_{jl}I_{2l}(t)}+\gamma R_j(t)-u_{3j}(t)S_j(t);\\
dW_j(t)/dt=kS_j(t) \sum_l{a_{jl}W_l(t)}-\beta^W_1 W_j(t) \sum_l{b_{jl}I_{1l}(t)}-\\
\hskip 38pt\beta^W_2 S_j(t) \sum_l{b_{jl}I_{2l}(t)}+u_{3j}(t)S_j(t)-\sigma_3W_j(t);\\
dI_{1j}(t)/dt= \beta^S_1 S_j(t) \sum_l{b_{jl}I_{1l}(t)}+\beta^W_1W_j(t)\sum_l{b_{jl}I_{1l}(t)}- \\
\hskip 38pt \varepsilon I_{1j}(t)\sum_l{b_{jl}I_{2l}(t)}-\sigma_1I_{1j}(t) -u_{1j}(t)I_{1j}(t);\\
dI_{2j}(t)/dt= \beta^S_2 S_j(t) \sum_l{b_{jl}I_{2l}(t)}+\beta^W_2W_j(t) \sum_l{b_{jl}I_{2l}(t)}+\\
\hskip 38pt \varepsilon I_{1j}(t)\sum_l{b_{jl}I_{2l}(t)}-\sigma_2I_{2j}(t) -u_{2j}(t)I_{2j}(t);\\
dR_j(t)/dt=\sigma_1I_{1j}(t)+u_{1j}(t) I_{1j}(t)+\sigma_2I_{2j}(t)+\\ 
\hskip 38pt u_{2j}(t)I_{2j}(t)+\sigma_3W_j(t)-\gamma R_j(t),
\end{array}
\end{equation}

}
\noindent where $\sum_l$ defines the sum from $1$ to $M$. Initial states are $S_j(0)>0,$ $W_j(0)\ge 0,$ $I_{1j}(0)>0,$ $I_{2j}(0)>0$, $R_j(0)=1-S_j(0)-W_j(0)-I_{1j}(0)-I_{2j}(0)$ for all clusters $j$.% $j\in \{1,\ldots,M\}.$

The aggregated system costs on the time interval $[0,T]$ are defined as $J=J_1+J_2$,
%\begin{equation}\label{functional_J}\begin{array}{c}
%J=J_1+J_2, \end{array}
%\end{equation}
where
{\small \begin{equation}\label{functional_SIR1}\begin{array}{c}
J_1=\int_0^T h_3\Big(\sum_j (u_{3j}(t))\Big)-L\Big(\sum_j W_j(t)\Big) dt,\\
J_2=\int_0^T \sum^2_{q=1} \bigg(f_q\big(\sum_j(I_{qj}(t))\big)+h_q\big(\sum_j (I_{qj}(t))\big)\bigg)-\\
\hskip 100pt g(\sum_j R_j(t)) dt.
\end{array}
\end{equation}}
and the optimal control problem is to minimize these costs, i.e., $\min_{\{u_{1j}, u_{2j}, u_{3j}\}} J.$

We  focus  on  a   case when both malware can cause extreme  damages, and there is a need to  lock down the entire system to prevent future destruction. To avoid this lockdown  or  other  expensive  security  activity, we have to construct a constant control such  that  any malware  will  be  instantly  eliminated,   even though the time when the viruses attack the system cannot be precisely identified. We assume that
\[ \max \left( h_1(u_{1j}), h_2(u_{2j}), h_3(u_{3j}), L(W_j), g(R_j) \right) << \]
\[ \min \left( f_1(I_{1j}), f_2(I_{2j})) \right), \] \[ \forall j, u_{1j}, u_{2j}, u_{3j}, W_j, R_j, I_{1j}, I_{2j}>0. \]

We have to define the condition for $u$ which remains system in disease free state with minimum costs. We assume that $h_1(u_{1j})=h_2(u_{2j})=h_3(u_{3j})=u$. The initial state of the system is the equilibrium point $E_2$ from the Section 3.2. (\ref{SIR2}) can be reformulated as:
{\small
\begin{equation}\label{cond_1}\begin{array}{l}
\beta^S_q S^0_j \sum_l{b_{jl}I_{ql}^0}+\beta^W_qW^0_j\sum_l{b_{jl}I_{ql}^0}+(-1)^q\varepsilon I_{1j}^0\sum_l{b_{jl}I_{ql}^0}- \\
\hskip 60pt \sigma_qI_{qj}^0 -u_{qj}(0)I_{qj}^0 \leq 0,\ q\in \{1,2\}.
\end{array}\end{equation}}
 It is assumed that  viruses can infect only one node at one time moment, then the system can be transformed in the following way:
{\small \begin{equation}\label{cond_2}\begin{array}{l}
\beta^S_q S^0_j b_{jm}+\beta^W_qW^0_j b_{jm}-  \sigma_q -u_{qj}(t) \leq 0,\ q=1,2,
\end{array}\end{equation}}
where $m$ is a node which was infected by a virus. Inequalities (\ref{cond_2}) can be rewritten as
{\small \begin{equation}\label{cond_3}\begin{array}{l}
u_{qj}(0) \geq (\beta^S_q S^0_j+\beta^W_qW^0_j) b_{jm}-  \sigma_1,\ q=1,2.\\
\end{array}\end{equation}}
We find control strategies that maintain the disease free state in the the worst case of epidemics. This value provides an estimation on system costs when $h_j=u_j(t)$ on the time interval $[0,T]$. Summing the control parameters gives:
{\small
\begin{equation}\begin{array}{l}
\hskip 52 pt \sum_j (u_{1j}(0)+u_{2j}(0)) \geq \\
\hskip 15pt \left((\beta^S_1+\beta^S_2 ) \sum_j S^0_j+(\beta^W_1+\beta^W_2 ) \sum_j W^0_j \right) b_{jm}-\\
\hskip 115pt M(\sigma_1+\sigma_2)=U,
\end{array} \end{equation}}

where $u_{ij}(t)$ is the control of a type $i \in \{1,2,3\}$ in a cluster $\mu$ at time $t$. As a result, we obtain
{\small
\begin{equation}\begin{array}{l}
J \rightarrow T \cdot \Big( \min (h_1(U),h_2(U)) - L(W(0)) - g(R(0)) \Big). \end{array} \end{equation}}

\vspace{-3mm}\section{Numerical Experiments}\label{numerical}
In this section, we present numerical case studies to  corroborate our results. For the experiments, we  use the following costs functions: infection costs -- $f_1(I_1(t))=30I_1(t)$ and $f_2(I_2(t))=40I_2(t)$; treatment costs -- $h_1(u_1(t))=20u^2_1(t)$, $h_2(u_2(t))=25u^2_1(t)$; vaccination cost -- $h_3(u_3(t))=10u^2_3(t)$; and utility functions are $L(W(t))=2W(t)$ and $g(R(t))=5R(t)$. The time interval in the first two experiments is equal to  [0,20].
\begin{figure}[h!]
\begin{center}
 \includegraphics[width=87mm]{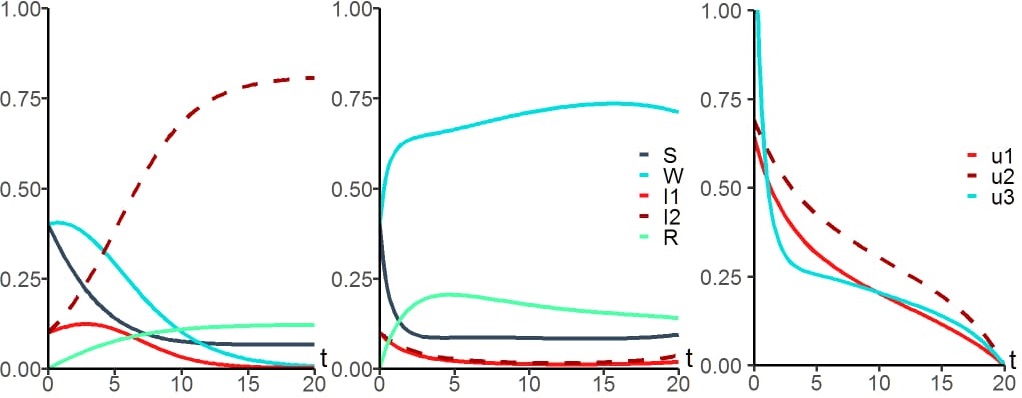}
\caption{Experiment I: Behavior of the system in the uncontrolled case (left), the controlled case (middle) and the structure of the optimal control(right). Parameters are: $k=0.3$, $\beta^S_1=0.35$, $\beta^S_2=0.45$, $\beta^W_1=0.25$, $\beta^W_2=0.35$, $\sigma_1 = 0.05$, $\sigma_2 = 0.03$, $\sigma_3 = 0.01$, $\gamma=0.2$, $\varepsilon = 0.5$).}
\end{center}\label{Exp1}
\end{figure}

Experiment I shows the behaviour of the SWIRS-model in two different cases: controlled and uncontrolled ones (Fig. 5). In the uncontrolled cases, at  $T=20$  the majority of nodes are infected by  virus $V_2$ ($I_2(20)=0.77$). The values of the functionals are equal to $J_1=-2.86$ and $J_2=10.41$. After the treatment and information dissemination about possible epidemic outbreaks, all infected nodes are cured. Here, all nodes are in the disease free state ($S(20)=0.29$, $W(20)=0.23$, $R(20)=0.48$) and values of the functionals are equal to $J_1=-11.12$ and $J_2=0.58$. Comparing the aggregated costs in the uncontrolled case ($J_{uncntl}=7.55$) and  the controlled case ($J_{cntl}=-10.54$), we can see that information spreading and applied treatment are beneficial.

Fig. 6 represents the dependence of the total number of infected nodes $I_{total}$ throughout the epidemic process on the parameters $k$ and $\sigma_3$ in the uncontrolled (left) and the controlled (right) cases, where $I_{total}=\int_0^T I_1(t)+I_2(t)dt.$
\begin{figure}[h!]
\begin{center}
 \includegraphics[width=80mm]{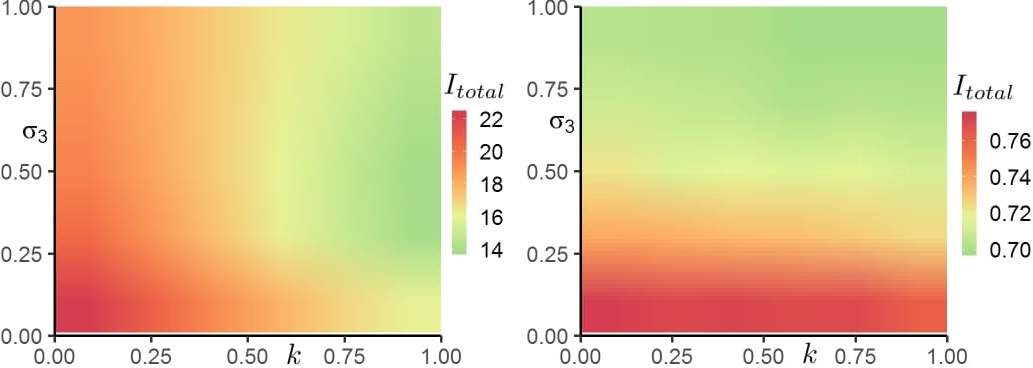}
\caption{Dependence of the total number of infected nodes on  parameters $k$ and $\sigma_3$.}
\end{center}\label{fig-color}
\end{figure}

In experiment II, we present the structure of the optimal control policies for the SWIR-model, when $\gamma=0$. In this case, after  the treatment, the recovered node will not be infected again during the contacts with  infected nodes. The final state of the system is ($0,0,0,0.82,0.18$). The aggregated system costs in the uncontrolled case are $J_{unctrl}=8.63$ ($J_1=-3.16$ and $J_2=11.79$). In the controlled case, the aggregated system costs reduced to $J_{ctrl}=-12.89$ ($J_1=-10.64$ and $J_2=-2.25$).
\begin{figure}[h!]
\begin{center}
  \includegraphics[width=87mm]{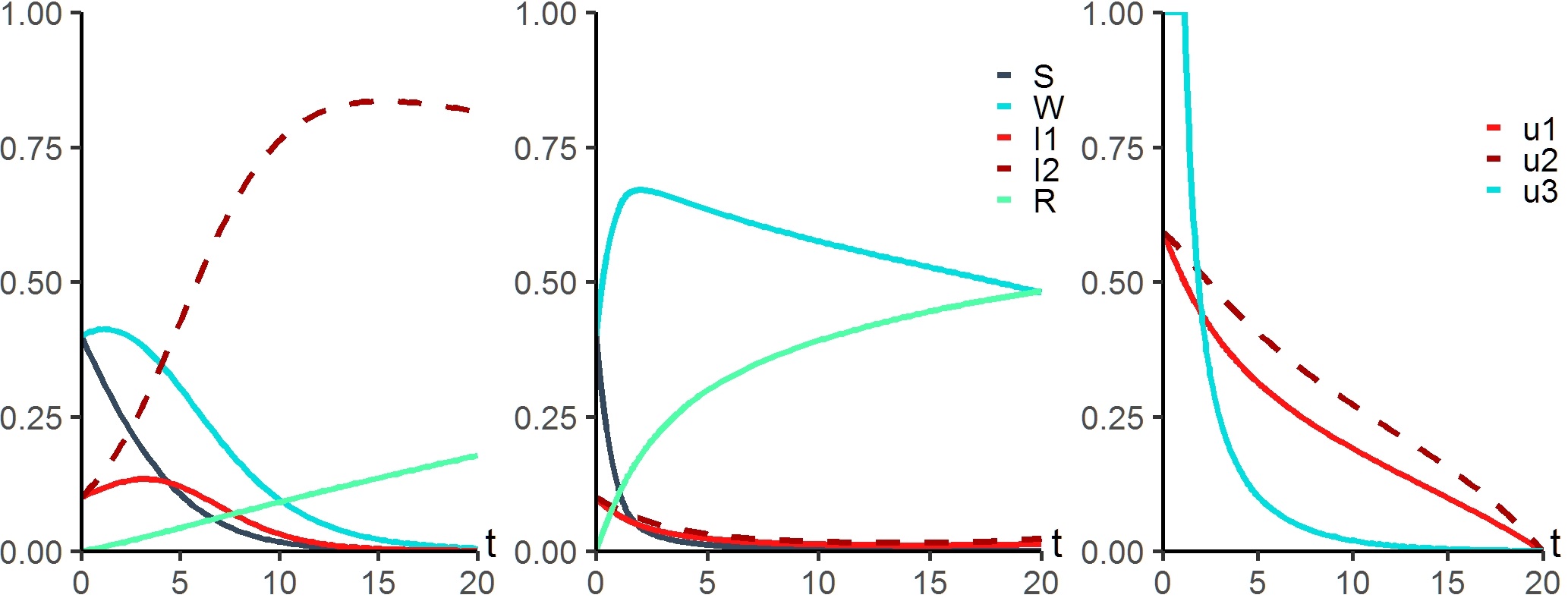}
  \caption{Experiment II: The behavior of the system in the uncontrolled case (left), the controlled case (middle), and the structure of the optimal control(right). Parameters are: $k=0.3$, $\beta^S_1=0.35$, $\beta^S_2=0.45$, $\beta^W_1=0.25$, $\beta^W_2=0.35$, $\sigma_1 = 0.05$, $\sigma_2 = 0.03$, $\sigma_3 = 0.05$, $\gamma=0$, $\varepsilon = 0.5$).}
\end{center}\label{Exp2}
\end{figure}

Experiment III presents the SWIRS model on a meta-population network, the behavior of the system (25) in two different clusters is represented in Fig. 8.   
\begin{equation}
A=B=
\begin{pmatrix}
1 & 0 \\
1 & 1
\end{pmatrix}.%\
\end{equation}
\begin{figure}[h!]
\begin{center}
  \includegraphics[width=70mm]{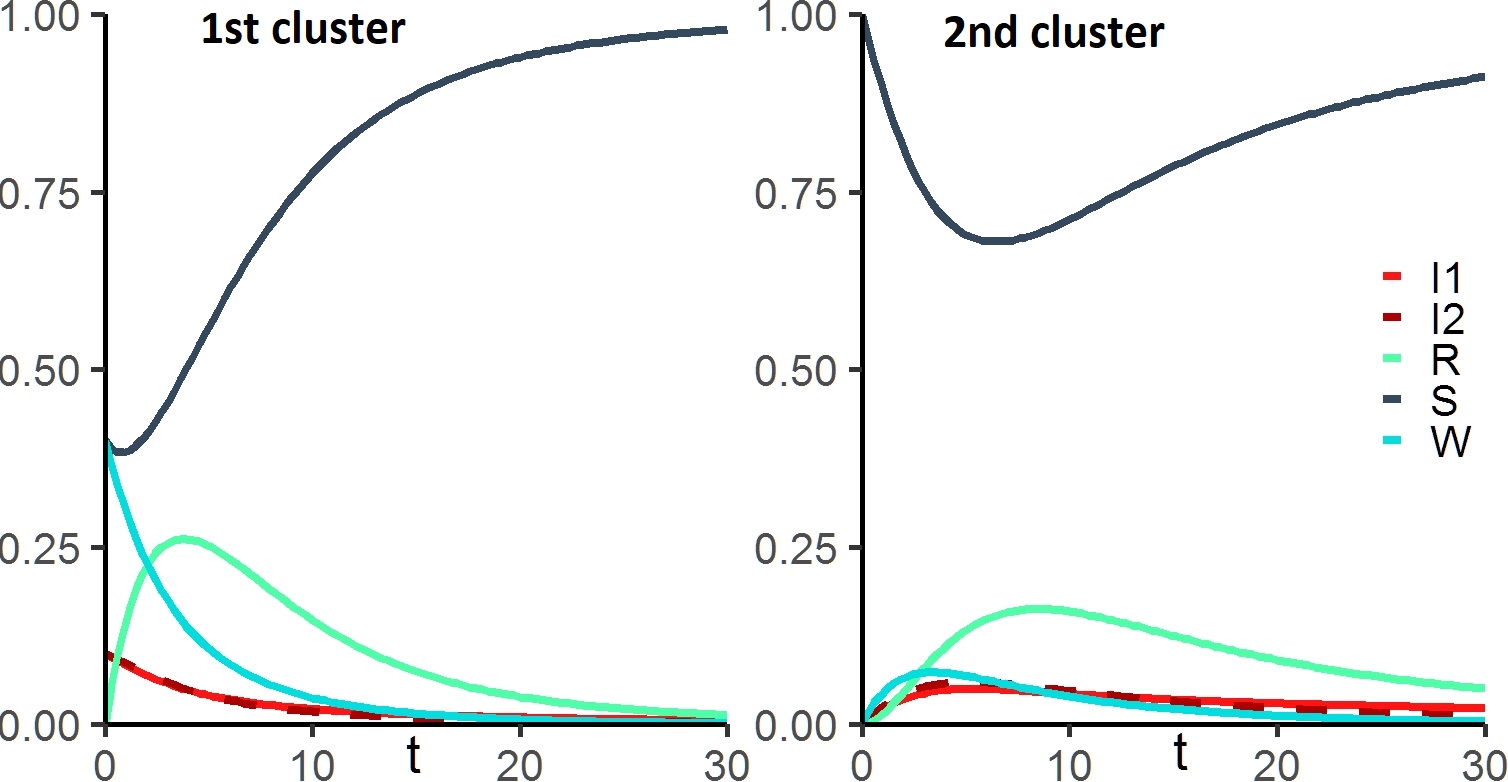}
  \caption{Experiment III: Behavior of the system in two different clusters of the population. Parameters are: $k=0.15$, $\beta^S_1=0.25$, $\beta^S_2=0.3$, $\beta^W_1=0.2$, $\beta^W_2=0.25$, $\sigma_1 = 0.3$, $\sigma_2 = 0.4$, $\sigma_3 = 0.3$, $\gamma=0.3$, $\varepsilon = 0.5$).}
\end{center}\label{fig:3}
\end{figure}

Initial parameters  are $X_1(0)=(0.4,0.4,0.1,0.1,0)$ and $X_2(0)=(1,0,0,0,0)$. Matrices A=B show the strong connections between these clusters, hence the epidemics which has been started in the first cluster continues in the second one. Final states are $X_1(30)=(0.97,0,0,0,0.03)$ and $X_2(30)=(0.9,0,0.02,0.02,0.06)$.

\section{Conclusions}\label{conclusion}
This paper presents a modified Susceptible-Warned- $\\$Infected-Recovered-Susceptible (SWIRS) model of simultaneous spreading of the virus protection information and the malware over a large population of nodes. We have investigated the stability of SWIR and SWIRS epidemic models with two coexisting malware types for heterogeneous populations. We have obtained the structure of the optimal control as well as the properties of feasible controls for a special class of cost functions. Numerical examples have been used to corroborate the results. We would further explore the stability properties of the epidemic process under optimal control. Another future work includes the extension of the SWIR model to an epidemic model over complex networks with different topologies.

\section{Acknowledgement}
The research has been partially supported by the RSF grant No. 16-19-10609, U.S. National Science Foundation Awards ECCS-1847056, CNS-1544782, and SES-1541164, and grant W911NF-19-1-0041 from U.S. Army Research Office (ARO).

\end{document}